\documentclass[11pt, a4paper]{article}
\usepackage{cite}
 \usepackage[utf8]{inputenc}
\usepackage{amssymb}
\usepackage{amsmath}
\usepackage{latexsym}
\usepackage{tikz}
\usepackage{sidecap,epstopdf}
\usepackage{subfigure}
\usepackage{color}
\usepackage{enumerate}

\usepackage{amsthm}\theoremstyle{remark}

\newcommand{\bte}{\begin{quote}\begin{theorem}}
\newcommand{\ete}[1]{\label{#1}\end{theorem}\end{quote}}
\newcommand{\bcom}{\begin{quote}\end{quote}}
\newcommand{\bex}{\begin{quote}\begin{example}}
\newcommand{\eex}[1]{\label{#1}\end{example}\end{quote}}
\newcommand{\bcon}{\begin{quote}\begin{conclusion}}
\newcommand{\econ}[1]{\label{#1}\end{conclusion}\end{quote}}
\newcommand{\bdefi}{\begin{quote}\begin{definition}}
\newcommand{\edefi}[1]{\label{#1}\end{definition}\end{quote}}

\newcommand{\blem}{\begin{quote}\begin{lemma}}
\newcommand{\elem}[1]{\label{#1}\end{lemma}\end{quote}}

\newcommand{\bpr}{\begin{quote}\begin{problem}}
\newcommand{\epr}[1]{\label{#1}\end{problem}\end{quote}}

\newcommand{\fr}{\frac}

\newcommand{\prt}{\partial}

\newcommand{\beq}{\begin{eqnarray}}
\newcommand{\eeq}[1]{\label{#1}\end{eqnarray}}
\newcommand{\bequ}{\begin{equation}}
\newcommand{\eequ}[1]{\label{#1}\end{equation}}
\newcommand\eq[1]{(\ref{#1})}
\newcommand{\bfi}{\begin{figure}[24]}
\newcommand{\efi}[1]{\caption{\label{#1}}\end{figure}}

\newcommand{\bfm}[1]{\mbox{\boldmath ${#1}$}}

\newcommand{\Bn}{{\textbf n}}

\newcommand{\Bx}{{\textbf x}}
\newcommand{\By}{{\textbf y}}
\newcommand{\Bz}{{\textbf z}}

\newcommand{\BC}{{\textbf C}}
\newcommand{\BD}{{\textbf D}}

\newcommand{\BI}{{\textbf I}}

\newcommand{\BO}{{\textbf O}}

\newcommand{\BS}{{\textbf S}}

\newcommand{\CF}{{\cal F}}
\newcommand{\CG}{{\cal G}}

\newcommand{\BGn}{\bfm\eta}

\newcommand{\BGx}{\bfm\xi}

\newcommand{\Ga}{\alpha}
\newcommand{\Gb}{\beta}
\newcommand{\Gd}{\delta}

\newcommand{\Gve}{\varepsilon}

\newcommand{\Gg}{\gamma}

\newcommand{\Gl}{\lambda}

\newcommand{\Go}{\omega}

\newcommand{\Gz}{\zeta}
\newcommand{\GD}{\Delta}

\newcommand{\GG}{\Gamma}

\newcommand{\GO}{\Omega}

\newcommand{\bexe}{\begin{quote}\begin{exercise}\inh}
\newcommand{\eexe}[1]{\label{#1}\end{exercise}\end{quote}}

\usepackage{graphics,graphicx}
\usepackage{color}

\begin{document}

\title{On meso-scale approximations for vibrations of  membranes with  lower-dimensional clusters of inertial inclusions}

\author{V.G. Maz'ya$^{1,2,3}$, A.B. Movchan$^{2}$, M.J. Nieves$^{4, 5}$  \\ ~
\\
$^1$ {\small Department of Mathematics, Link\"oping University,} \\ {\small Link\"oping S--581 83, Sweden}  \\
$^2$ {\small Department of Mathematical Sciences, University of Liverpool, } \\
{\small Liverpool L69 7ZL, UK} \\
$^3$ {\small RUDN University, 6 Miklukho-Maklay St, Moscow, 117198, Russia} \\
$^4$ {\small School of Computing and Mathematics, Keele University,} \\
{\small Staffordshire, ST5 5BG, UK}\\
$^5$ {\small  Department of Mechanical, Chemical and Material Engineering,} \\ {\small University of Cagliari, Cagliari, 09123, Italy}}

\date{}

\maketitle

\begin{abstract}
In this paper we consider formal asymptotic algorithms for a class of meso-scale approximations for problems of vibration of elastic membranes, which contain clusters of small inertial inclusions distributed along contours of pre-defined smooth shapes.  
Effective transmission conditions have been identified for inertial structured interfaces, and approximations to  solutions of eigenvalue problems have been derived for domains containing lower-dimensional clusters of inclusions.
\end{abstract}

\vspace{0.5in}

\centerline{\it In honour of Professor N.N. Uraltseva}

\vspace{0.5in}

\section{Introduction}


We address a class of asymptotic approximations for models of vibrations of two-dimensional elastic membranes, containing  clusters of small inclusions. These clusters are assumed to be distributed along one-dimensional sets. 
The approach is based on the method of meso-scale asymptotic approximations \cite{MM, MMN, MMN1, MMN2}, and the  inertia of inclusions has been taken into account. 
The method of meso-scale asymptotic approximations was first introduced in \cite{MM}  for asymptotic problems in domains with large clusters of small inclusions, and it provides an efficient alternative to homogenisation approximations (see, for example, \cite{MKhr}), especially for the cases where clusters are inhomogeneous, non-periodic, and where the size of inclusions is in the meso-scale range compared to the distance between the inclusions. The asymptotic approximations  obtained in \cite{MM} involve a linear combination of solutions to certain model problems whose coefficients satisfy a linear algebraic system. The solvability of this system was proved under weak geometrical assumptions, and both uniform and energy estimates for the remainder term were derived.
 
Meso-scale approximations for eigenvalue problems in domains with clusters of many inclusions were analysed in \cite{Meso_eigen}. Fundamental ideas of the method of compound asymptotic expansions \cite{MNP} in domains with singularly perturbed boundaries were used.

Although problems of wave scattering in the low frequency regime in solids with many inclusions can be addressed by the method of homogenisation, important features linked to the wave scattering from individual inclusions may require the use of the dynamic Green's functions. In particular, Green's function for the Helmholtz equation in a periodic domain was analysed in \cite{Linton, Kurkcuetal}. A formal  procedure was introduced in \cite{Foldy} for isotropic scattering by randomly distributed scatterers separated by a finite distance. Analysis of waves in a plane with  semi-infinite arrays of isotropic scatterers is discussed in  \cite{LinMart}. Waves governed by the Helmholtz equation in a doubly periodic array with an elementary cell containing several scatterers were analysed in \cite{SchCr}, based on the approach of  \cite{Foldy}, \cite{LinMart} and asymptotics representing singular perturbation leading-order approximations, similar to simplified cases of  \cite{MMN}, \cite{Meso_eigen} and \cite{MNP}. Dispersion of waves analysed in \cite{SchCr}  shows dynamic anisotropy linked to the scatterers within the elementary cell.

\begin{figure}
\centering
\includegraphics[width=0.6\textwidth,keepaspectratio]{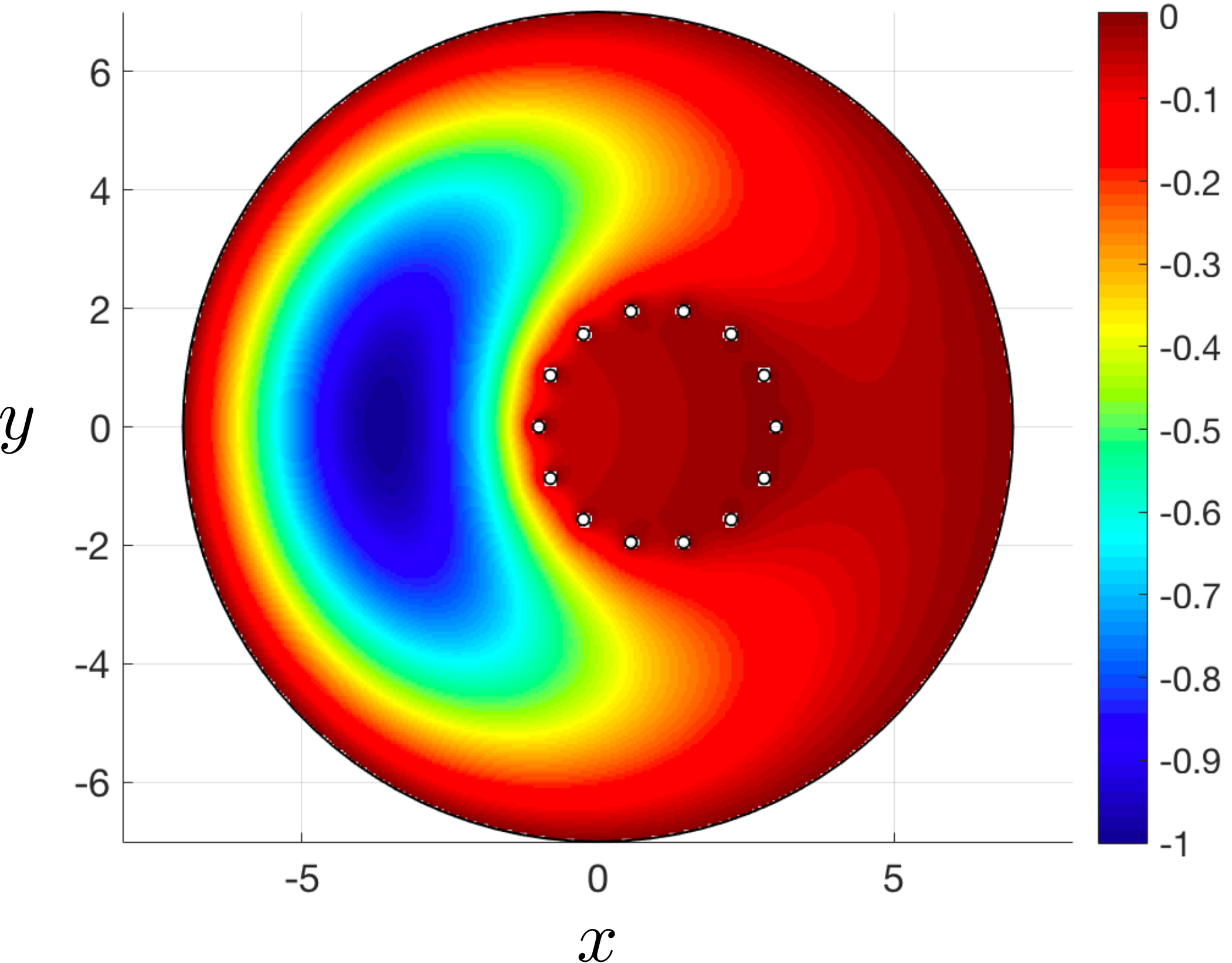}
\caption{Cluster of small inclusions placed along a contour. The figure shows the eigenfunction corresponding to the first eigenvalue of the Dirichlet problem for the Laplacian, the eigenvalue is 
$0.30816$. The main region is the disk of radius $R=7$ with the centre at the origin. Small circular rigid inclusions of radii $r=0.1$ have their centres along the circle of radius $2$ with the centre at $(1, 0)$. The computation is produced in COMSOL.}
\label{fig1}
\end{figure}

For example, when a cluster of $N$ small circular inclusions of radius $\Gve$ is assumed to be placed along a simple smooth closed curve in the plane we can use the method of meso-scale asymptotic approximations for  solutions of scattering problems.   
When an incident wave of a radian frequency $\Go$ is generated by a remote time-harmonic point source placed at $\By \in \mathbb{R}^2$ the approximation of the total field is constructed in the form
$$
\mathcal{G}(\Bx, \By, \omega) \sim G(|\Bx- \By|, \omega)+\sum_{j=1}^N \beta_j(\By) U^{(j)}(\Bx)\;,
$$
where $G(|\Bx- \By|, \omega)$ is Green's function for the Helmholtz equation, $U^{(j)}$ are special model fields associated with individual small inertial inclusions placed at $\BO^{(j)}$, and $\Gb_{j}$ are the coefficients, found by solving an  algebraic system of equations. 
The proof of the solvability of this algebraic system is similar to \cite{MM}.
Compared to meso-scale two-dimensional clusters, this is a lower--dimensional case where inclusions are distributed along a contour. 
In addition to the above point-wise asymptotic approximation, the approach discussed here, can also be used to derive an effective transmission problem for the domain containing an inertial structured interface, and its solution provides a homogenisation approximation, which takes into account inertia of small inclusions.  In particular, when small inclusions of mass $m$ form a periodic cluster of the overall mass $M$ distributed over a circle $\GG$ of unit radius, the coefficients $\Gb_j$ can be approximated by 
$$
\beta_j(\By)=\mathfrak{G}(\BO^{(j)}, \By)\;, \quad 1\le j \le N\;,
$$
with the function $\mathfrak{G}(\BO^{(j)}, \By)$ being the solution of the  transmission problem for the following equation 
$$
\mu \GD_x \mathfrak{G}(\Bx, \By)  + \rho \Go^2 \mathfrak{G}(\Bx, \By)  + {\frak R}(m, M, \Go)
\Gd(\Bx - \By) = 0,
$$
where $\mu$ is the stiffness coefficient, $\rho$ is the mass density, and ${\frak R}(m, M, \Go)$ is the coefficient, that depends on the inertial properties of the cluster. While $\mathfrak{G}(\BO^{(j)}, \By)$ satisfies the radiation condition at infinity, the transmission conditions across $\GG$ have the form
$$
\big[\mathfrak{G}(\Bx, \By)\big]=0\;,\quad \Big[\frac{\partial \mathfrak{G}}{\partial r}(\Bx, \By)\Big] =-\frac{M }{2 \pi m} \Ga_\Gve 
{\frak R}(m, M, \Go) \mathfrak{G}(\Bx, \By)\;. 
$$
Here $\Ga_\Gve = -4{\rm i \mu}/H^{(1)}_0  (\sqrt{\rho/\mu} \Go {\varepsilon}),$ and $H^{(1)}_0$ is the Hankel function of the first kind. 


For non-resonance forced problems, we also discuss solutions which can be interpreted as fields associated with the time-harmonic scattering from a lower-dimensional cluster of small inertial inclusions. The incident field is produced by a localised source, associated with a forced motion of one of rigid inclusions within the array, and  formally the algorithm requires an asymptotic approximation of the harmonic capacitary potential in the low frequency regime.  The approximation corresponds to a non-resonance case. 

An example of a cluster of inclusions placed along a curve in a two-dimensional elastic membrane is shown in Fig. \ref{fig1}. Time-harmonic vibrations are considered here, and the inertia of small inclusions, as well as their size, separation and the size of the cluster represent parameters of the multi-scale asymptotic approximation. 

We also analyse asymptotic solutions of an eigenvalue problem for a domain containing a lower-dimensional cluster of small inertial inclusions. 
Computations, discussed in the text, give a comparison between an analytical asymptotic procedure and a FEM simulation produced in COMSOL Multi-Physics shown in Fig. \ref{fig1}, with the analytical logarithmic asymptotic  approximation shown in Fig. \ref{fig2} and produce a remarkably good result. 

The structure of the paper is as follows. In Section \ref{1Dcluster} we introduce the asymptotic model of a one-dimensional cluster of many inertial inclusions placed along a closed simple contour, which forms an inertial structured interface. Section \ref{Qstatic} includes analysis of a low-frequency scattering in a finite elastic membrane with a sparse cluster of small movable rigid inclusions. The approximation of the first eigenvalue and the corresponding eigenfunction for a sparse or one-dimensional cluster of inertial inclusions in a finite elastic membrane is discussed in Section \ref{eigen_cluster}.

\section{Elastic membrane with a one-dimensional cluster of inclusions}   
\label{1Dcluster}

\vspace{0.1in}The problem considered concerns the scattering  phenomena due to a cluster of $N$ inertial inclusions, $F_\varepsilon^{(j)}$ with masses $m_j$  and centres $\BO^{(j)}$, $1\le j \le N$, in a membrane that has density $\rho$ and shear modulus $\mu$.  Here $\Gve$ is a small positive parameter, characterising the relative diameter of $F_\varepsilon^{(j)}$, similar to \cite{MMN}. The overall mass of the cluster $M=\sum_{j=1}^Nm_j$\;, is finite.
Here we seek the asymptotic approximation of the solution of the scattering problem for a wave initiated by a point source positioned outside the cluster
\begin{equation}\label{eqGeps1}
\mu\Delta \mathcal{G}(\Bx, \By, \omega)+\rho \omega^2 \mathcal{G}(\Bx, \By, \omega)+\delta(\Bx-\By)=0\;, \quad \Bx \in \mathbb{R}^2 \backslash \cup_{j=1}^N {\overline{F}^{(j)}_\varepsilon}\;,
\end{equation}
\begin{equation}\label{eqGeps2}
\mathcal{G}(\Bx, \By, \omega)=C_j\;, \quad \Bx \in \partial F^{(j)}_\varepsilon\;,
\end{equation}
\begin{equation}\label{eqGeps3}
\frac{\partial \mathcal{G}}{\partial r_x}(\Bx, \By, \omega)-{\rm i}k_0 \mathcal{G}(\Bx, \By, \omega)=O\Big(\frac{1}{r_x^{3/2}}\Big)\;, \quad \text{ as }\quad  r_x\to \infty \;,\end{equation}
where the position of the point force is given by ${\bf y}\in \mathbb{R}^2\backslash \cup_{j=1}^N {\overline{F}^{(j)}_\varepsilon}$, with $\By$ being separated by a finite distance from the cluster itself,  $r_x=|\Bx|$ and $C_j$ are constants that are determined in what follows. Here, $\omega$ and $k_0=\sqrt{{\rho \omega^2}/{\mu}}$ denote the radian frequency and wavenumber, respectively, of waves in the medium. In addition, as each mass vibrates in the membrane, the linear momentum balance for the mass is taken into account, imposing an additional condition in the form
\begin{equation}\label{eqlinbalance}
-m_j \omega^2 C_j =\mu \int_{\partial F_\varepsilon^{(j)}}\frac{\partial \CG}{\partial n}(\Bx, \By, \omega) ds_{\Bx}\;,
\end{equation}
for $1\le j \le N$, and $\Bn$ being the unit outward normal with respect to $F_\varepsilon^{(j)}$. 
This problem can be interpreted as that of the time-harmonic scattering of a wave produced by a point source in the presence of a cluster of small inertial inclusions positioned along a a simple smooth curve in a plane.

\subsection{Model problems}\label{modprob}

We use the algorithm of the method of meso-scale asymptotic approximations \cite{MMN} to 
formally construct  an asymptotic approximation of the field $\mathcal{G}(\Bx, \By, \omega)$.
The procedure requires several model problems. 

\subsubsection{Dynamic Green's function}
By $G$ we denote the dynamic Green's function for the infinite membrane, that satisfies the equation
\begin{equation}
\mu \Delta G(|\Bx-\By|, \omega)+\rho \omega^2 G(|\Bx- \By|, \omega )+\delta(\Bx-\By)=0\;,\qquad  \Bx, \By \in \mathbb{R}^2\;, \label{Grep0}
\end{equation}
and the outgoing wave solution has the representation 
\begin{equation}\label{Grep}
G(|\Bx -\By|, \omega)=\frac{{\rm i}}{4\mu} H^{(1)}_0\Big(k_0 |\Bx-\By|\Big)\;,
\end{equation}
{where $k_0=\sqrt{{\rho \omega^2}/{\mu}}$, and $H^{(1)}_0$ 
is the Hankel function of the first kind}. In particular, 
\[
G(|\Bx-\By|, \omega) \sim -\fr{1}{2 \pi \mu} \log (k_0 |\Bx-\By| ),
\]
as $k_0 |\Bx-\By| \to 0.$

\subsubsection{Dirichlet problem for the Helmholtz equation in the exterior of a finite inclusion}
To allow for the correction of discrepancies on interior boundaries we introduce the problem
\begin{equation} \label{U0}
\mu\Delta U^{(j)}(\Bx)+\rho\omega^2 U^{(j)}(\Bx)=0\;, \quad \Bx \in \mathbb{R}^2 \backslash {\overline{F}_\varepsilon^{(j)}}\;, 
\end{equation} 
\begin{equation}\label{U1}
U^{(j)}(\Bx)=1\;, \quad \Bx \in \partial F_\varepsilon^{(j)}\;,
\end{equation} \label{U1a}
\hspace{-0.1cm}where the $U^{(j)}$ also satisfies the radiation condition at infinity:
\begin{equation}
\frac{\partial U^{(j)}}{\partial r}(\Bx)-{\rm i}k_0 U^{(j)}(\Bx)=O(\frac{1}{r^{3/2}})\;, \quad \text{ as }\quad  r\to \infty \;.\end{equation}

\subsection{Meso-scale approximation - formal asymptotics} 
The leading order approximation for the  solution of (\ref{eqGeps1})--(\ref{eqlinbalance}) is sought in the form
\begin{equation}
\label{CGapp}
\mathcal{G}(\Bx, \By, \omega) \sim G(|\Bx- \By|, \omega)+\sum_{j=1}^N \beta_j(\By) U^{(j)}(\Bx)\;,
\end{equation}
where the coefficients $\beta_j$ are to be determined.

\subsubsection{The algebraic system}

It follows from   the boundary conditions  (\ref{eqGeps2}) that to leading order  we have
\begin{equation}\label{BC1}
C_k={G}(|\BO^{(k)}- \By|, \omega)+\beta_k(\By)+\sum_{\substack{j\ne k\\1\le j \le N}} \beta_j(\By)U^{(j)}({\BO^{(k)}})\;,\quad \Bx \in \partial F^{(k)}_\varepsilon\;, 
\end{equation}
with $1\le k \le N$.
Using the ``equations of motions'' (\ref{eqlinbalance}) for  individual inclusions together with (\ref{CGapp}), we deduce that 
the constants $C_k$, $1\le k\le N$, also satisfy the relations
\bequ
-\fr{m_k  \Go^2}{\mu}  C_k = \int_{\prt F_\Gve^{(k)}} \fr{\prt \CG(\Bx, \By, \Go)}{\prt n} d s_x \nonumber
\end{equation}
\begin{eqnarray}
&=& \int_{\prt F_\Gve^{(k)}} \fr{\prt G(|\Bx-\By|, \Go)}{\prt n} d s_x + \Gb_k(\By) 
 \int_{\prt F_\Gve^{(k)}} \fr{\prt U^{(k)}(\Bx)}{\prt n} d s_x  \nonumber \\
&+& \sum_{\substack{j \ne k\\ 1\le j \le N}} \beta_j (\By)   \int_{\prt F_\Gve^{(k)}} \frac{\partial U^{(j)} (\Bx) }{\partial n}d s_{x}, \label{ED1}
\end{eqnarray}
where $\Bn$ stands for the unit outward normal on $\prt F_\Gve^{(j)}.$ Taking into account \eq{Grep0}, \eq{U0}, \eq{U1}, and when $j \neq k$, using the integration by parts one can represent the boundary integrals in \eq{ED1} via the integrals over $F_\Gve^{(j)}$, and hence the representation \eq{ED1} takes the form
\beq
{m_k  \Go^2} C_k &=&  -\mu \Gb_k(\By) 
 \int_{\prt F_\Gve^{(k)}} \fr{\prt U^{(k)}(\Bx )}{\prt n} d s_x   \label{ED2} \\
&+& {\rho \Go^2} \int_{F_\Gve^{(k)}} \Big\{  G(|\Bx-\By|, \Go) +\sum_{\substack{j \ne k\\ 1\le j \le N}} \beta_j (\By)   U^{(j)} (\Bx) \Big\} d \Bx . \nonumber 
\end{eqnarray}
Combining \eq{BC1} and \eq{ED2}, we derive the algebraic system of equations for the coefficients $\Gb_j, ~ 1 \leq j \leq N, $ in the asymptotic approximation \eq{CGapp}
$$
m_k \Go^2 G(|\BO^{(k)}- \By|, \omega) - {\rho \Go^2} \int_{F_\Gve^{(k)}}  G(|\Bx-\By|, \Go) d \Bx
$$
\bequ
+ \Gb_k(\By) \Big\{  m_k \Go^2 + \mu   \int_{\prt F_\Gve^{(k)}} \fr{\prt U^{(k)}(\Bx )}{\prt n} d s_x \Big\}
\eequ{ED3}
$$
+ m_k \Go^2 \sum_{\substack{j\ne k\\1\le j \le N}} \beta_j(\By) \{  U^{(j)}({\BO^{(k)}})
- \fr{\rho}{m_k} \int_{F_\Gve^{(k)}} U^{(j)} (\Bx) d \Bx \} =0.
$$

The justification of solvability of the above algebraic system for the coefficients $\Gb_j$ is similar to \cite{MM}, and we do not discuss it here.

In particular, when inclusions $F_\Gve^{(j)}$ are circular with radii $\varepsilon r^{(j)}$ and  centres $\BO^{(j)}$, where $\varepsilon$ is a small non-dimensional parameter, we have 
\begin{equation}\label{Uj}
U^{(j)}(\Bx)=\alpha^{(j)}_\varepsilon G(|\Bx - \BO^{(j)}|, \omega)\;,
\end{equation}
where 
\[\alpha^{(j)}_\varepsilon=-\frac{4{\rm i \mu}}{H^{(1)}_0\Big(k_0{\varepsilon}r^{(j)}\Big)} \;.\]

In this case, the algebraic system \eq{ED3} for the coefficients $\beta_j$, $1\le j \le N$ becomes
\begin{eqnarray}
&& m_k  {G}(|\BO^{(k)}- \By|, \omega)-{{{\rho}}\int_{ F^{(k)}_\varepsilon }{ G}(|\Bx- \By|, \omega)d {\Bx}} \nonumber \\
&&+\beta_k\Big(m_k +\frac{2\pi  {\mu } \varepsilon r^{(k)} \Ga_\Gve^{(k)}}{\omega^2}\frac{\partial G(r, \Go)}{\partial r}\Big|_{r=\varepsilon r^{(k)}}\Big)\nonumber\\
&& +\sum_{\substack{j \ne k\\ 1\le j \le N}}\beta_j \Ga_\Gve^{(j)}\Big\{m_k G(|\BO^{(k)}-\BO^{(j)}|,\Go) \nonumber \\
&&-{ {\rho}\int_{F^{(k)}_\varepsilon }  { G}(|\Bx-\BO^{(j)}|,\Go) d {\Bx}}\Big\}=0\;.\nonumber \\\label{algeqs}
\end{eqnarray}

 \subsection{Scattering by a cluster of identical inclusions placed along a circular contour}

Here, we assume that the circular inclusions $F_\varepsilon^{(j)}$ have the same radii $\Gve$ and mass $m$ and that their centres are uniformly  distributed along the unit circular contour $\Gamma_1 = \{ \Bx: |\Bx| = 1\}$. Also, let $d = 2 \pi  /N$, and assume that $\Gve < \mbox{Const} ~ d^{3/2}.$ 
We analyse the case when $N \to \infty$, while $M = m N$ is fixed, and obtain 
a problem concerning a membrane with an inertial ring.

When the inclusions $F_\varepsilon^{(k)}$, $1\le k \le N$, are all circular, some of the integrals appearing in (\ref{algeqs}) can be evaluated explicitly. In particular, when $\By$ is separated by a finite distance from the cluster of the inclusions, 
using the  Graf's addition theorem and expanding the Hankel function in a series we have:
\[H_0^{(1)}(k_0|\Bx-\By|)= \sum_{\nu=-\infty}^\infty H_\nu^{(1)}(k_0|\By-\BO^{(k)}|) J_\nu (k_0 |\Bx - \BO^{(k)}|)e^{{\rm i}\nu (\pi -\theta_{k}+\theta_{ k, \By})}\;.\]
where $\theta_{k}$  is the polar angle of $\Bx$ measured with respect to $\BO^{(k)}$ and $\theta_{k, \By}$ is the polar angle of $\BO^{(k)}$ measured with respect to $\By$. Hence, as $\Gve \to 0,$
\begin{eqnarray*}
\int_{F_\varepsilon^{(k)}} G(|\Bx- \By|, \omega) ds_{\Bx}&=& \pi \Gve^2 G( |\BO^{(k)} - \By|, \omega) + O(\Gve^3) = O(\Gve^2),
\end{eqnarray*}
and for $j \neq k$, it is derived that 
\begin{eqnarray*}
\int_{F_\varepsilon^{(k)}} G(|\Bx -\BO^{(j)}|, \omega) ds_{\Bx}&=& \pi \Gve^2 G( |\BO^{(k)} - \BO^{(j)}|, \omega) 
\end{eqnarray*}
$$
+ O(\Gve^3/d) = O(\Gve^2 |\log d|).
$$
We also note that 
\[ 
{2\pi \varepsilon \mu r^{(k)}}\frac{\partial G}{\partial r}(r, \omega)\Big|_{r =
\varepsilon r^{(k)}} =-1 + O(\Gve^2 |\log \Gve| ).\]
For the  cluster of identical inclusions, the system  (\ref{algeqs}), to leading order, becomes
\begin{eqnarray}
&& m { G}( |\BO^{(k)} - \By|, \omega) +\beta_k (\By) \Big(m - \frac{\Ga_\Gve}{\omega^2}
\Big)\nonumber\\
&&+ m \Ga_\Gve \sum_{\substack{j \ne k\\ 1\le j \le N}}\beta_j (\By) G(|\BO^{(k)}-\BO^{(j)}|)=0\;, \nonumber \\\label{algeqs1nc}
\end{eqnarray}
where $1 \leq k \leq N,$ and
\[\alpha_\varepsilon=-\frac{4{\rm i \mu}}{H^{(1)}_0\Big(k_0{\varepsilon}\Big)} \;.\] 

\vspace{.2in}

%
%
\subsection{Derivation of the transmission conditions for $\mathfrak{G}$ in the auxiliary problem}

Assume that positions of the centres of inclusions $F^{(j)}_\varepsilon$ are   
\[\BO^{(j)}=(\cos(2\pi (j-1)/N ), \sin(2\pi (j-1)/N ) )^{\rm T}\;, \qquad  1\le j\le N\;,\]
%
and let $\mathfrak{G} (\Bx, \By)$ be a function such that 
\begin{equation}\label{criterion}
\beta_j(\By)=\mathfrak{G}(\BO^{(j)}, \By)\;, \quad 1\le j \le N\;.
\end{equation}
%
%
%
%
%
Equations \eq{algeqs1nc} can be re-written in the form
\begin{eqnarray}
&& 
m { G}( |\BO^{(k)} - \By|, \omega) +\beta_k (\By) \Big(m - \frac{\Ga_\Gve}{\omega^2}
\Big)
\nonumber\\
&&
+ \fr{M \Ga_\Gve }{2 \pi}\sum_{\substack{j \ne k\\ 1\le j \le N}}\beta_j (\By) G(|\BO^{(k)}-\BO^{(j)}|) \fr{2 \pi}{N}=0\;.\label{algeqs1ncA}
\end{eqnarray}
By considering the sum in \eq{algeqs1ncA} as the Riemann sum and taking the limit as $N \to \infty$, we arrive at
%
\[
m \Go^2 { G}( |\Bx - \By|, \omega) + \mathfrak{G}(\Bx, \By) (m \Go^2 - \Ga_\Gve) 
\]
\bequ
+\frac{M\omega^2 \Ga_\Gve}{2\pi}\int_0^{2\pi} G(|\Bx- \BGn|, \Go)\mathfrak{G}(\BGn, \By) d\theta_{\scriptsize{\BGn}}=0\;. 
\eequ{algeqs1ncB}
When this equation is extended to $\Bx$ outside $\GG_1$, one can apply the Laplacian in $\Bx$, use \eq{Grep0} and derive 
\[
0 = -m \Go^2 (\fr{\rho \Go^2}{\mu}G (|\Bx - \By|, \Go) + \fr{1}{\mu} \Gd(\Bx-\By)) +  (m \Go^2 - \Ga_\Gve) \GD_x \mathfrak{G}(\Bx, \By)
\]
\bequ
- \fr{M \Go^4 \Ga_\Gve \rho}{2 \pi \mu}  \int_{\GG_1}  G(|\Bx- \BGn|, \Go)\mathfrak{G}(\BGn, \By) d\theta_{\scriptsize{\BGn}}.
\eequ{eqfrakGa}
Furthermore, using  \eq{algeqs1ncB} and \eq{eqfrakGa} we obtain the equation for $\mathfrak{G}(\Bx, \By)$ in the form
\bequ
\mu \GD_x \mathfrak{G}(\Bx, \By)  + \rho \Go^2 \mathfrak{G}(\Bx, \By)  + \fr{m \Go^2}{\Ga_\Gve - m \Go^2}
\Gd(\Bx - \By) = 0. \label{eqfrakG}
\end{equation}
We note that $|\Ga_\Gve| = O(|\log \Gve|^{-1})$, with $m$ being of order $O(d)$, and the radian frequency $\Go^2$ serving as an additional control parameter.


In addition, we take the normal derivative in (\ref{algeqs1ncB}) and use the following relations when 
$\Bx \in \Gamma_1$
\[
\lim_{\scriptsize{\BGx} \to \Gamma_1^\pm}\int_{\Gamma_1} \frac{\partial G}{\partial r}(|\BGx- \BGn|, \Go)\mathfrak{G}(\BGn, \By) d s_{\BGn}=\mp \frac{1}{2}\mathfrak{G}(\Bx, \By)+\int_{\Gamma_1} \frac{\partial G}{\partial r}(|\Bx- \BGn|, \Go)\mathfrak{G}(\BGn, \By) ds_{\scriptsize{\BGn}}\;, \quad 
\]
(see \cite{SV}), and obtain the transmission condition across $\GG_1$
\bequ
\big[\mathfrak{G}(\Bx, \By)\big]=0\;,\quad \Big[\frac{\partial \mathfrak{G}}{\partial r}(\Bx, \By)\Big] = 
\frac{M \omega^2 \Ga_\Gve}{2 \pi (m \Go^2 - \Ga_\Gve)}\mathfrak{G}(\Bx, \By)\;, 
\eequ{jump_c}
with $\big[\mathfrak{G}(\Bx, \By)\big]$ denoting the jump of $\mathfrak{G}(\Bx, \By)$ across $\Gamma_1$, while the point $\By$ is separated by a finite distance from $\GG_1$:
\[\big[\mathfrak{G}(\Bx, \By)\big]=\mathfrak{G}(\Bx, \By)\Big|_{\Bx \in \Gamma_1^+}-\mathfrak{G}(\Bx, \By)\Big|_{\Bx \in\Gamma_1^-}\;.\]
Thus, the function $\mathfrak{G}(\Bx, \By)$ can be defined as a solution of \eq{eqfrakG}, \eq{jump_c}, subject to the radiation condition at infinity. This provides an alternative homogenisation approximation 
\eq{criterion} for the coefficients $\Gb_j$ in \eq{CGapp}, which takes into account the inertial transmission conditions across the structured interface formed by a cluster of small inclusions. 

\section{Low-frequency scattering in a finite elastic membrane with a sparse cluster of small movable 
rigid inclusions}
\label{Qstatic}

Consider a finite elastic membrane $\Omega \subset \mathbb{R}^2$, together with a finite sparse cluster  of small rigid inclusions $F_\varepsilon^{(j)} \subset \GO, ~j = 1,\ldots,N,$ of zero mass density, containing interior points $\BO^{(j)}$. It is assumed that the small inclusions $F_\varepsilon^{(j)}, j = 1,\ldots,N,$ are separated by the finite distance from the exterior boundary $\prt \GO.$
We use the notation $$\GO_N = \Omega \setminus \cup_{j=1}^N \overline{F}_\Gve^{(j)}.$$
Given the mass density $\rho$ and  a time-harmonic vibration of small radian frequency $ \omega,$ applied to the inclusion $F_\varepsilon^{(1)}$, the amplitude of the out-of-plane displacement satisfies the problem: 

\begin{eqnarray}
\mu \GD u(\Bx) + \rho \Go^2 u(\Bx) &=& 0, ~~ \Bx \in \GO_N
\label{E1} \\
u(\Bx) &=& 1, ~~  \Bx \in \partial F_\varepsilon^{(1)}  \label{E2}  \\
u(\Bx) &=& A_j,  ~ \Bx \in \partial F_\varepsilon^{(j)} ~ 1 < j \leq N,   \label{E3} \\
u(\Bx) &=& 0, ~~  \Bx \in \partial \GO \label{E4} 
\end{eqnarray}
where 
\bequ
\int_{\prt  F_\varepsilon^{(j)}} \fr{\prt u}{\prt n} d s = 0, ~~ 1 < j \leq N.
\eequ{E5}
The constants $A_j$ are to be determined, and it is assumed  that small inclusions $F_\varepsilon^{(j)}$  are separated by a finite distance. We introduce a small non-dimensional parameter $f$, with $\Gve < f$, in such a way that
$\rho \Go^2 = f \Gl$. Here, we use the normalisation $\mu=1$, and hence the equation \eq{E1} takes the form
\bequ
\GD u(\Bx) + f \Gl u(\Bx) = 0, ~~ \Bx \in \GO_N.
\eequ{E6}
Also, the notation $F^{(j)}$ is used for scaled inclusions, such that $F^{(j)} = \{\Bx: \Gve \Bx + \BO^{(j)} \in F^{(j)}_\Gve\}.$

\subsection{Green's function and the relative  capacitary potential}

We use the result of \cite{MMN} and employ the {\em relative capacitary potential} $P_\Gve^{(1)}(\Bx)$, which satisfies the boundary value problem in $\GO\setminus \overline{F}_\Gve^{(1)}$
\begin{eqnarray}
\GD P_\Gve^{(1)}(\Bx) &=& 0, ~~ \Bx \in \GO\setminus \overline{F}_\Gve^{(1)},
\label{E7} \\
P_\Gve^{(1)}(\Bx) &=& 1, ~~  \Bx \in \partial F_\varepsilon^{(1)},  \label{E8}  \\
P_\Gve^{(1)}(\Bx) &=& 0, ~~  \Bx \in \partial \GO .\label{E9} 
\end{eqnarray}

The notations $G(\Bx, \By)$ and $g(\Bx, \By)$ are used for Green's functions in $\GO$ and $\mathbb{R}^2 \setminus \overline{F}^{(1)}$, respectively.  Thus, 
\begin{eqnarray}
\GD G(\Bx, \By) + \Gd(\Bx - \By) &=& 0, ~~ \Bx, \By \in \GO,
\label{E10} \\
G(\Bx, \By) &=& 0, ~~  \Bx \in \partial \GO, ~\By \in \GO, \label{E11} 
\end{eqnarray}
and 
\begin{eqnarray}
\GD g(\BGx, \BGn) + \Gd(\BGx - \BGn) &=& 0, ~~ \BGx, \BGn \in \mathbb{R}^2 \setminus \overline{F}^{(1)},
\label{E12} \\
g(\BGx, \BGn) &=& 0, ~~  \BGx \in \partial F^{(1)},  \BGn \in \mathbb{R}^2 \setminus \overline{F}^{(1)}, \label{E13} \\
g(\BGx, \BGn) ~\mbox{is bounded} & \mbox{as} & |\BGx| \to \infty ~\mbox{and}~ \BGn \in \mathbb{R}^2 \setminus \overline{F}^{(1)}. \label{E14}
\end{eqnarray}

According to Lemma 1.4 in \cite{MMN}, the relative capacitary potential  $P_\Gve^{(1)}$ has the asymptotic representation
\bequ
P_\Gve^{(1)}(\Bx) = \fr{-G(\Bx, \BO^{(1)}) + \Gz(\fr{\Bx-\BO^{(1)}}{\Gve})-\fr{1}{2 \pi} \log \fr{|\Bx - \BO^{(1)}|}{\Gve r_F}}{\fr{1}{2 \pi } \log \fr{\Gve r_F}{R_\GO}} + p_\Gve(\Bx),
\eequ{E16}
where $p_\Gve(\Bx) = O(\Gve |\log \Gve|^{-1})$ uniformly with respect to $\Bx \in \GO \setminus \overline{F}_\Gve^{(1)}$. In the above formula \eq{E16}, the notations $r_F$ and $R_\GO$ stand for the inner conformal radius of $F^{(1)}$ with respect to $\BO^{(1)}$, and the outer conformal radius of $\GO$ relative to $\BO^{(1)}$, respectively, as discussed in Section 1.2.1 of \cite{MMN}. The function $\Gz$ is defined by
\bequ
\Gz(\BGn) = \lim_{|\BGx| \to \infty} g(\BGx, \BGn).
\eequ{E17}
In particular, for the case when the inclusion $F_\Gve^{(1)}$ is circular of radius $\Gve r_F$, the formula \eq{E16} simplifies as
\bequ
P_\Gve^{(1)}(\Bx) = {\frak B}_1 G(\Bx, \BO^{(1)})  + O(\Gve), ~~ {\frak B}_1 = -(\fr{1}{2 \pi } \log {\Gve r_F} + H(\BO^{(1)}, \BO^{(1)}))^{-1},
\eequ{E16a}
where the regular part $H(\Bx, \By)$ 
of Green's function $G$ is a harmonic function defined by 
\bequ
H(\Bx, \By) = \fr{1}{2 \pi} \log |\Bx - \By|^{-1} - G(\Bx, \By).
\eequ{E15}

\subsection{Formal asymptotic approximation}

The asymptotic approximation is sought in the form
\bequ
u(\Bx)\sim P^{(1)}_\Gve (\Bx) + f u^{(2)}(\Bx),
\eequ{E18}
and, owing to \eq{E3},  the constants $A_j, ~ j = 2,\ldots, N,$ representing the rigid motion displacements of massless inclusions, are
\bequ
A_j \sim P_\Gve^{(1)} (\BO^{(j)}) + f A_j^{(2)},
\eequ{E19} 
whereas $u^{(2)}$ satisfies the boundary value problem
\beq
\GD u^{(2)}(\Bx) + \Gl P^{(1)}_\Gve (\Bx) &=& 0 ~~ \mbox{in} ~ \GO_N, \label{E20} \\
u^{(2)} (\Bx) &=& 0, ~ \Bx \in \prt F_\Gve^{(1)}, \label{E21} \\
u^{(2)}(\Bx) &=& A^{(2)}_j, ~ \Bx \in \prt F_\Gve^{(j)},  j= 2, \ldots, N, \label{E22} \\
u^{(2)}(\Bx) &=& 0, ~ \Bx \in \prt \GO. \label{E23}
\end{eqnarray}
By assuming a circular shape of the inclusion $F_\Gve^{(1)}$, using formula \eq{E16a}, and introducing an auxiliary problem
\beq
\GD V(\Bx) + \Gl {\frak B}_1 G (\Bx, \BO^{(1)}) &=& 0 ~~ \mbox{in} ~ \GO, \label{E24} \\
V(\Bx) &=& 0, ~ \Bx \in \prt \GO, \label{E25}
\end{eqnarray}
we deduce
\bequ
V(\Bx) = \Gl {\frak B}_1 \int_{\GO} G(\Bz, \Bx) G(\Bz, \BO^{(1)}) d \Bz,
\eequ{E26}
and
\bequ
A^{(2)}_j = V(\BO^{(j)}) - V(\BO^{(1)}) P_\Gve^{(1)}(\BO^{(j)}). 
\eequ{E27}
In this case, the approximation \eq{E18} can be rewritten in the form
\bequ
u(\Bx) = P^{(1)}_\Gve (\Bx) + f ( V(\Bx) - V(\BO^{(1)}) P^{(1)}_\Gve (\Bx) ) + {\frak r}(\Bx),
\eequ{E28}
where the remainder term satisfies the problem
\beq
\GD {\frak r}(\Bx) + f \Gl {\frak r}(\Bx) &=& -f \Gl (P_\Gve^{(1)} (\Bx) - {\frak B}_1 G(\Bx, \BO^{(1)})) \label{E28a} \\
&-& f^2 \Gl (V(\Bx) - V(\BO^{(1)}) P_\Gve^{(1)} (\Bx) ), \nonumber \\
{\frak r}(\Bx) &=& 0 ~~ \mbox{on} ~ \prt \GO,  \label{E28b} \\
{\frak r}(\Bx) &=& -f (V(\Bx) - V(\BO^{(1)}))  ~~ \mbox{on} ~ \prt F_\Gve^{(1)},  \label{E28c} \\
{\frak r}(\Bx) &=& A_j - P_\Gve^{(1)}(\Bx) - f (V(\Bx) - V(\BO^{(1)}) P_\Gve^{(1)}(\Bx))   \nonumber \\
&& ~~ \mbox{on} ~ \prt F_\Gve^{(j)}, ~ j = 2, \ldots, N.  \label{E28d} 
\end{eqnarray}
For the case when the inclusion $F_\Gve^{(1)}$ is circular, taking into account that the small quantity $f \Gl$ is separated from the spectrum, and using the formulae
 \eq{E19}, \eq{E27} and \eq{E28} we deduce
that  
 \bequ
 {\frak r}(\Bx) = O(f  \Gve)
 \eequ{E29}
when $\Bx$ is outside a neighbourhood of the cluster of small inclusions, and
 \bequ
 {\frak r}(\Bx) = O({\rm max} \{f  \Gve |\log \Gve|, ~\Gve\})
 \eequ{E29a}
when $\Bx$ is in the vicinity of the cluster of small inclusions.

The idea of the proof for the case of a sparse cluster of small circular inclusions is based on the representation of the solution $ {\frak r}$ of \eq{E28a}--\eq{E28d} as a sum of integrals over $\GO_N$ and the boundaries of  $F_\Gve^{(j)}$ :
$$ 
{\frak r}(\Bx) = - \int_{\GO_N} \CF(\By) g_\GO(\Bx, \By) d \By
$$
$$+ \sum_{j=1}^N
\int_{\prt F_\Gve^{(j)}} \Big\{
g_\GO(\Bx, \By) \fr{\prt {\frak r}}{\prt n_\By}(\By) - \Phi_j(\By) \fr{\prt g_\GO}{\prt n_\By}(\Bx, \By)
\Big\} d s_\By 
$$
where $\CF$ is the right-hand side in \eq{E28a}, $\Bn$ is the unit outward normal with respect to $\GO_N$, and the functions  $\Phi_j$ are the right-hand sides in \eq{E28c} for $j=1$, and in \eq{E28d} for $ 2 \leq j \leq N.$ Here $g_\GO(\Bx, \By)$ is Green's function for the Helmholtz equation in the unperturbed domain $\GO$:
$$
\GD_\Bx g_\GO(\Bx, \By) + k^2 g_\GO(\Bx, \By) + \Gd(\Bx - \By) = 0, ~ \Bx, \By \in \GO, 
$$
$$
g_\GO(\Bx, \By) = 0, ~\mbox{when} ~ \Bx \in \prt \GO, ~\By \in \GO,
$$
where $k = \sqrt{f \Gl}.$
We also note that the flux of ${\frak r}$ over the boundaries of small inclusions is not zero, i.e.
$$
\int_{F_\Gve^{(1)}} \fr{\prt {\frak r}}{\prt n} d \Bx = O(f), ~\int_{F_\Gve^{(j)}} \fr{\prt {\frak r}}{\prt n} d \Bx = O(f \Gve^2), ~j=2,\ldots, N.
$$
By considering three cases for ${\frak r}(\Bx)$, (a) $\Bx$ is outside the neighbourhood  of the  sparse  cluster of small circular inclusions, (b)  $\Bx$ is in the vicinity of $F_\Gve^{(1)}$, and (c) 
$\Bx$ is in the vicinity of $F_\Gve^{(j)}, ~ j= 2,\ldots,N$, we deduce \eq{E29} and \eq{E29a}.

\section{The first positive eigenvalue for a lower-- dimensional cluster of inertial inclusions in a finite elastic membrane}
\label{eigen_cluster}
Section \ref{Qstatic} has addressed  low frequency non-resonance vibrations for a sparse cluster of massless inclusions (each individual inclusion has zero inertia). In this case, Green's function for the Laplacian  is used as one of the model solutions. 
If the small inclusions have a non-zero inertia the problem requires Green's function for the Helmholtz operator, which takes into account time-harmonic vibrations of the multi-scale system. 
In this section, we consider the eigenvalue problem for a finite membrane $\Omega \subset \mathbb{R}^2$, and contained in this membrane is a sparse or lower-dimensional cluster of small inertial inclusions $F_\varepsilon^{(j)}$ of the mass $m_j$, $1\le j \le N$ .  In this case, if $d$ is the minimum distance between two neighbouring inclusions within the cluster, we assume that $dN = O(1).$ As before, we use the notation $\Omega_N:= \Omega \backslash \cup_{j=1}^N {\overline{F}^{(j)}_\varepsilon}$. 
The  eigenfunction $u_N$ and the corresponding eigenvalue $\lambda_N$ are defined as the solution of the  problem
\begin{equation}\label{fineig1}
\Delta u_N(\Bx)+\lambda_N u_N(\Bx)=0\;, \quad \Bx \in \Omega_N\;,
\end{equation}
\begin{equation}\label{fineig2}
u_N(\Bx)=0\;, \quad \Bx \in \partial \Omega ,
\end{equation}
\begin{equation}\label{fineig3}
u_N(\Bx)=A_j\;, \quad \Bx \in \partial F^{(j)}_\varepsilon\;,
\end{equation}
with 
\begin{equation}\label{fineig3a}
-\gamma_j\lambda_N A_j=\int_{\partial F_\varepsilon^{(j)}}\frac{\partial u_N(\Bx)}{\partial n}ds\;, \quad \Bx \in \partial F^{(j)}_\varepsilon\;.
\end{equation}
where $\lambda_N= \rho \omega^2/\mu$ and $\gamma_j=\frac{m_j}{\rho}$, $1\le j \le N$, and $n$ is the unit outward normal with respect to $F^{(j)}_\varepsilon$.
The quantities $\Gg_j$ represent the equivalent inertial area, which would be required for the inclusion of mass density $\rho$ to have the mass $m_j$. For convenience, we also introduce the normalisation $\mu=1$ and the notation $\Gl = \rho \Go^2,$ where $\Go$ is the radian frequency of the time-harmonic vibrations and $\rho$ is the mass density of the membrane. It is also assumed that $\sum_{j=1}^N m_j = M = O(1),$ and thus $\gamma_j = O(M/(N \rho))$.


Here, we are interested in evaluating the first positive eigenvalue, which takes into account the inertia of the small inclusions $F_\Gve^{(j)}.$ 


\subsection{Model problem: Green's function for a finite membrane}
Assuming a non-resonance regime for $\Gl$, we shall use Green's function $G_\Omega(\Bx, \By, {\lambda})$ in $\GO$, i.e. the solution of the following problem
\begin{equation}\label{fineig4}
\Delta G_\Omega(\Bx, \By, {\lambda})+\lambda G_\Omega(\Bx, \By, {\lambda})+\delta(\Bx-\By)=0\;, \quad \Bx \in \Omega\;,
\end{equation}
\begin{equation}\label{fineig5}
G_\Omega
(\Bx, \By, {\lambda})=0\;, \quad \Bx \in \partial \Omega\;,
\end{equation}
and the regular part $R_\Omega$ of $G$ is defined as
\begin{equation}\label{fineig6}
R_\Omega(\Bx, \By, {\lambda})=-\frac{{1}}{4}Y_0(\sqrt{\lambda}|\Bx-\By|)-G_\Omega(\Bx, \By, {\lambda})\;,
\end{equation}
where $Y_0 (\sqrt{\Gl} r)$ is the Bessel function of the second kind, and
$$
Y_0 (\sqrt{\Gl} r) \sim  \fr{2}{ \pi} \log (\sqrt{\Gl} r) , ~~\mbox{as} ~ \sqrt{\Gl} r \to 0.
$$
Here, we also use the auxiliary functions $V^{(j)}$ (compare with \eq{U0},\eq{U1})  defined by 
\begin{equation} \label{UU0}
\Delta V^{(j)}(\Bx, \Gl)+\Gl V^{(j)}(\Bx, \Gl)=0\;, \quad \Bx \in \GO \backslash {F_\varepsilon^{(j)}}\;, 
\end{equation}  
\begin{equation}\label{UU1}
V^{(j)}(\Bx, \Gl)=1\;, \quad \Bx \in \partial F_\varepsilon^{(j)}\;,
\end{equation}
\begin{equation}\label{UU2}
V^{(j)}(\Bx, \Gl)=0\;, \quad \Bx \in \partial \GO\;.
\end{equation}
When $F_\Gve^{(j)}$ are small circular inclusions of radii $\Gve r_j$, $V^{(j)}$  are approximated in the form
\bequ
V^{(j)}(\Bx, \Gl) = \Gb_\Gve^{(j)} G_\GO (\Bx, \BO^{(j)}, \Gl) + O(\Gve |\log \Gve|^{-1}),
\eequ{E30}
where
\bequ
\Gb_\Gve^{(j)} = -\Big( \frac{{1}}{4}Y_0(\sqrt{\lambda} \Gve r_j)+R_\Omega(\BO^{(j)}, \BO^{(j)}, {\lambda})     \Big)^{-1} . 
\eequ{E31}

\subsection{Formal approximation of the  first eigenvalue and corresponding eigenfunction}
We look for an approximation $\lambda_N^{(0)}$ of the first eigenvalue of \eq{fineig1}--\eq{fineig3a} and an approximation $\sigma_N(\Bx, \lambda_N^{(0)})$ of the  corresponding eigenfunction.
We use the representation
\begin{equation}\label{sigN}
\sigma_N(\Bx, \lambda_N^{(0)})=\sum_{j=1}^N C_j V^{(j)} (\Bx, \Gl_N^{(0)}), 
\end{equation}
with the coefficients $C_j$ being unknown intensities of inertia forces exerted by the inclusions on the membrane. The approximation \eq{sigN} satisfies \eq{fineig1}, with $\Gl$ being replaced by $\Gl_N^{(0)}$,  and the boundary condition  \eq{fineig2} on $\prt \GO$. 


The boundary conditions \eq{fineig3},  \eq{fineig3a} yield the system of algebraic equations for the coefficients  $C_j$.

\subsubsection{The algebraic system}

When $\Bx \in \prt F_\Gve^{(k)}$, we use  \eq{sigN}, together with \eq{fineig3},  \eq{fineig3a}, to deduce the following relations
\bequ
-\fr{1}{\Gg_k \Gl_N^{(0)}} T^{(k)}[ \sigma_N] = C_k + \sum_{\substack{j \neq k \\ 1 \leq j \leq N}} C_j  V^{(j)} (\Bx, \Gl_N^{(0)}), ~
\Bx \in \prt F_\Gve^{(k)}, ~ k = 1, \ldots, N,
\eequ{E32}
where
\bequ
T^{(k)}[ \sigma_N]  = \int_{\prt F_\Gve} \fr{\prt \sigma_N}{\prt n} ds.
\eequ{E33}

Taking into account that 
\bequ
T^{(k)} [ \sigma_N]  = -\Gb_\Gve^{(k)} C_k + O (\Gve |\log \Gve|),
\eequ{E34}
and expanding $V^{(j)} (\Bx, \Gl_N^{(0)})$ in the series form near $\Bx = \BO^{(k)}, ~ k \neq j$, we arrive at
the algebraic system equations for the coefficients $C_k$
\bequ
C_k (1 - \fr{\Gb_\Gve^{(k)}}{\Gg_k \Gl_N^{(0)}}) + \sum_{\substack{j \neq k \\ 1 \leq j \leq N}}  C_j V^{(j)} (\BO^{(k)}, \Gl_N^{(0)}) =0, ~ k = 1, \ldots, N,
\eequ{E35}
and furthermore, taking into account \eq{E30} we approximate the coefficients $C_j$ as non-trivial solutions of the following homogeneous system
\bequ
\big\{\BI - \BD_\Gve (\Gl_N^{(0)}) + \BS_\Gve (\Gl_N^{(0)}) \big\}  \BC = 0,
\eequ{E36}
where $\BC = (C_1, ... , C_N)^T,$ $\BI$ is the identity matrix, and $\BD_\Gve (\Gl_N^{(0)}),  \BS_\Gve (\Gl_N^{(0)})$ are $N \times N$ matrices defined as
\beq
\BD_\Gve (\Gl_N^{(0)}) &=& \fr{1}{\Gl_N^{(0)}}\mbox{diag} \big\{  \fr{\Gb_\Gve^{(1)}}{\Gg_1},...,\fr{\Gb_\Gve^{(N)}}{\Gg_N }  \big\},  \nonumber \\
~ (\BS_\Gve (\Gl_N^{(0)}))_{kj} &=& \Gb_\Gve^{(j)} G_\GO (\BO^{(k)}, \BO^{(j)}, \Gl_N^{(0)}) (1-\Gd_{jk}).
\label{E37}
 \end{eqnarray}
 The approximation $\Gl_N^{(0)}$ of the first eigenvalue is defined from the equation
 \bequ
 \mbox{det} \big\{\BI - \BD_\Gve (\Gl_N^{(0)}) + \BS_\Gve (\Gl_N^{(0)}) \big\} = 0.
 \eequ{E38}
 
 Given $\Gl_N^{(0)}$, the approximation $\sigma_N$ of the eigenfunction is defined by \eq{sigN}, with the coefficients $C_j$ obtained from  \eq{E36}. 
 
 We also note that the above algorithm equally applies to the Dirichlet problem \eq{fineig1}--\eq{fineig3}, with $A_j$ being zero in the right-hand side of \eq{fineig3}. In this case, the matrix term $\BD_\Gve$ in \eq{E36} is replaced by zero matrix. 
 

 \begin{figure}
\centering
\includegraphics[width=0.6\textwidth,keepaspectratio]{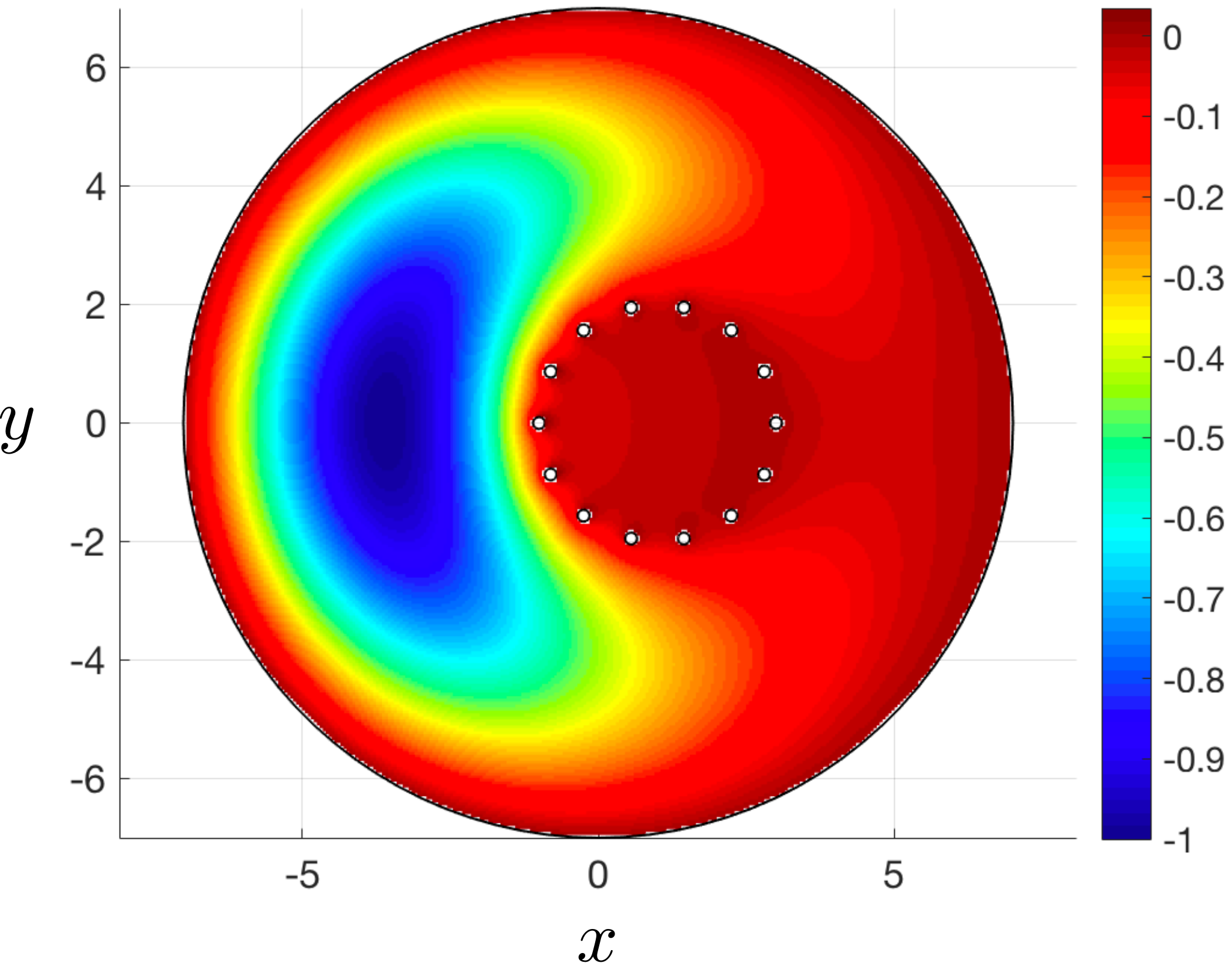}
\caption{Asymptotic approximation of the eigenfunction corresponding to the first eigenvalue of the Dirichlet problem in the domain with the  cluster of small inclusions placed along a contour. The approximation of the first eigenvalue is $\Gl_N^{(0)} = 0.30678$. The main region is a disk of radius $R=7$. Small circular rigid inclusions of radii $r=0.1$ have their centres along the circle of radius $2$, the same as in Fig. \ref{fig1}. }
\label{fig2}
\end{figure}
 
 \subsection{An example}
 
 The above scheme is applicable to sparse clusters or dense one-dimensional clusters of inertial inclusions.
 When the masses $m_j, ~j = 1,\ldots, N,$ of small inclusions increase, with the membrane mass density $\rho$ remaining constant, coefficients  $\Gg_j$ in \eq{fineig3a} also increase and the first eigenvalue of \eq{fineig1}--\eq{fineig3a} decreases accordingly.
 
 However, we obtain a different problem in the limit when the masses $m_j$ tend to infinity, with $\rho$ being finite. The corresponding formulation will be the Dirichlet eigenvalue problem \eq{fineig1}--\eq{fineig3}, with the homogeneous boundary conditions   \eq{fineig3}. The equation \eq{E36}, with $\BD_\Gve$ being zero, defines the approximation for the first eigenfrequency. 
 
 Here, we discuss an example, where the results of the asymptotic approximation are compared to  an independent numerical simulation in COMSOL Multi-Physics.
 
 When the domain $\GO$ is the disk $\{\Bx: |\Bx| < R\}$ of radius $R$ with the centre at the origin, and $\By \neq {\bf 0},$ the regular part of Green's function in \eq{fineig6} and \eq{E31} is evaluated with the use of Graf's addition formula, as follows
 \bequ
 R_\GO(\Bx, \By, \Gl) = \mbox{Re} \Big\{  \sum_{n=-\infty}^\infty \Ga_n J_n (\sqrt{\Gl} |\Bx|) e^{{\rm i} n  \theta_{\bf x,0}}  \Big\},
 \eequ{E39}
 where the coefficients $\Ga_n$ are given in the form
 \bequ
 \Ga_n = \fr{{\rm i} H_n^{(1)} (\sqrt{\Gl} R) J_n(\sqrt{\Gl} |\By|)}{4 J_n(\sqrt{\Gl} R)} e^{- {\rm i} n (\theta_{\bf 0,y} +\pi)},
  \eequ{E39a}
 where $\theta_{\bf x,0}$ is the polar angle of $\Bx$ with respect to the centre ${\bf 0}$ of the disk, and $\theta_{\bf 0,y}$ is the polar angle of ${\bf 0}$ with respect to $\By$.
 
 When $\By = {\bf 0}$, we have
 \bequ
  R_\GO(\Bx, {\bf 0}, \Gl) = - \fr{Y_0(\sqrt{\Gl} R)}{4 J_0 (\sqrt{\Gl} R)} J_0 (\sqrt{\Gl} |\Bx|).
 \eequ{E40}
 
 Figure \ref{fig2} shows the asymptotic approximation of the eigenfunction corresponding to the first eigenvalue of the Dirichlet problem in the domain $\GO_N$ with the  cluster of small inclusions placed along a contour. The approximation of the first eigenvalue is $
 0.30678$. An independent numerical simulation produced in COMSOL Multi-Physics gives the first eigenvalue of $\Gl_N = 0.30816$, with an excellent agreement observed in Fig. \ref{fig1} and Fig. \ref{fig2}.
 
 
 \section{Concluding remarks}

 We have given an outline of problems, where meso-scale asymptotic approximations provide an analytical insight to analysis of time-harmonic wave problems in two-dimensional domains with small inclusions arranged in a sparse cluster or a one-dimensional cluster placed along a curve.
 
We  note that logarithmic asymptotics, required for analysis of Dirichlet eigenvalue problems in domains with small inclusions often lead to a constraint of exponentially small size of the inclusions. On the other hand, if the rigid small inclusions are considered as ``movable'', the inertia of the inclusion is used as an additional control parameter in the meso-scale asymptotic approximation.

Three classes of formulations discussed here show several directions where asymptotics of  time-harmonic waves prove to be useful. In particular, in Section  \ref{1Dcluster}, for a finite mass cluster of small inertial inclusions we have derived a transmission condition across an inertial structured interface. While in Section \ref{Qstatic}, the quasi-static Green's function has been successfully employed, together with the relative capacitary potential, it is essential to take into account wave scattering and reflection and to engage Green's function for the Helmholtz operator in the analysis of the eigenvalue problem for a domain with a lower-dimensional cluster of small inclusions in Section \ref{eigen_cluster}. As illustrated, the derived asymptotic formulae provide a constructive analytical tool and are straightforward to use in practical examples and computations.

\section*{Acknowledgements}
V.G.M. acknowledges that this publication has been prepared with the support of the ``RUDN University Program 5-100''. A.B.M would like to thank the EPSRC (UK) for its support through the Programme Grant no. EP/L024926/1.  M.J.N  gratefully acknowledges the support of the EU H2020 grant MSCA-IF-2016-747334-CAT-FFLAP.

\end{document}